\documentclass{amsart}

\usepackage{graphicx}
\usepackage{amsfonts}
\usepackage{amssymb}
\usepackage{xcolor}

\newtheorem{theorem}{Theorem}[section]
\newtheorem{conjecture}[theorem]{Conjecture}

\newcommand{\pres}[3]{\textnormal{#1} \langle #2 \: | \: #3 \rangle}

\title[The B. B. Newman Spelling Theorem]
 {The B. B. Newman Spelling Theorem} 

\author{Carl-Fredrik Nyberg-Brodda*}
\address{School of Mathematics, University of East Anglia, Norwich NR4 7TJ, England, UK}
\email{c.nyberg-brodda@uea.ac.uk}
\thanks{*The author is currently a PhD student at the University of East Anglia, UK}

\begin{document}

\begin{abstract}
This article aims to be a self-contained account of the history of the B. B. Newman Spelling Theorem, including the historical context in which it arose. First, an account of B. B. Newman and how he came to prove his Spelling Theorem is given, together with a description of the author's efforts to track this information down. Following this, a high-level description of combinatorial group theory is given. This is then tied in with a description of the history of the word problem, a fundamental problem in the area. After a description of some of the theory of one-relator groups, an important part of combinatorial group theory, the natural division line into the torsion and torsion-free case for such groups is described. This culminates in a statement of and general discussion about the B. B. Newman Spelling Theorem and its importance.
\end{abstract}

\maketitle

Mathematical research is done by mathematicians. It is all too easy to forget that whereas the truth of mathematical statements is independent of the people who proved them, the process of research leading up to these proofs is very human. When reading the title of this essay, one may expect to find herein the statement of the B. B. Newman Spelling Theorem, and an explanation and history of some of the relevant terms; perhaps an exposition of the history of the mathematical context in which the theorem appears is also to be expected. All of this is present. However, one may also justifiably ask who B. B. Newman was. Beginning this article is therefore a section on the very human story of the mathematician behind the theorem, and how my efforts to track him down instead led to me being tracked down by him. It is not expected of the reader to understand all of the mathematics within. The heart of this story does not lie with its mathematical content.

\section{B. B. Newman}

I was first made aware of the B. B. Newman Spelling Theorem following a meeting with my doctoral supervisor at the very beginning of my PhD, and while it would be wholly incorrect to claim that this initial encounter resulted in any deep understanding of the theorem, it left me in no doubt as to the importance of the theorem. In 1968, when the theorem was first proven, the machinery used today to prove the theorem was not as well oiled as it is today. My curiosity was therefore piqued as to what the original proof may have looked like, na\"ively expecting a relatively straightforward path to finding it. It is difficult to overstate how wrong these expectations would turn out to be.

In any case, the optimism was not necessarily unjustified, as there was not one, but two natural starting points for my search. In most citations two references were cited for the Spelling Theorem: a bulletin article from 1968, and B. B. Newman's doctoral thesis from 1968, submitted to the University of Queensland (UQ), in Australia. While the bulletin article certainly contained the statement of the Spelling Theorem, and many of its consequences, it, understandably, did not contain the proof I sought. Neither was any light yet shone on the mysterious initials ``B. B.'', and so my attention moved towards finding the doctoral thesis. Writing to the university library at UQ, the amount of information available to me on this topic was significantly increased after some correspondence. Bill Bateup Newman, hereinafter simply referred to as Bill, was born in 1936, and graduated with, among other degrees, an M.Sc. from the School of Mathematics at UQ in 1964, and was on record for completing a PhD in 1970. The master's thesis, entitled ``Almost Just Metabelian Groups'', was freely available in the archives, and a copy was sent to me. However, there was no trace of any doctoral thesis. Furthermore, the mystery thickening, there was not even any sign as to who might have been his doctoral supervisor. After much effort from the librarian and archivists, little else was found in the archives. Slightly defeated, I decided to slump into my chair to at least try and figure out what an \textit{almost just metabelian} group might be.

As I was reading through the preface of the master's thesis, a line suddenly caught my eye, in which Bill thanked his supervisor, a Dr M. F. Newman. The trail was picked up again, and looking around the staff pages of various Australian universities, I was able to find Professor Emeritus Michael (Mike) Frederick Newman at the Australian National University. I wrote to him, almost, but not entirely, convinced that this surely must be the same mathematician, and a couple of days later received a response. ``I am the M. F. Newman you are looking for'' he began reassuringly. While Mike had not had any contact with Bill for many years, nor even knew whether he was still alive, nor knew the whereabouts of any copies of his thesis, he did inform me of something which took me by surprise. While the doctoral studies had been for a degree at UQ, the doctoral supervisor of Bill had been Gilbert Baumslag. This was surprising for a number of reasons. The first was that Baumslag, as will become clear from the number of mentions of the same in the following sections, was a very major name in the area, and a quick search revealed that there was no record\footnote{After discovering the connection, I submitted an entry to the Mathematics Genealogy Project (https://genealogy.math.ndsu.nodak.edu), which now lists Bill as a student of Baumslag's.} online of him ever having supervised Bill. The second was that Baumslag had been based at City College New York (CCNY), in a city almost antipodal to Queensland, for a large portion of his life. Of course, there was only one thing I could do at this point. I wrote to CCNY, and asked them whether they had a copy of Bill's thesis. A few days later, I received an email. 

It would not be an exaggeration to say that I have never been as surprised by an email as I was upon seeing who had sent it. There, in my inbox, was an email from B. B. Newman. It was beyond words; what had come to pass was that for all my efforts to track this mathematician down, he had managed to find me first. A myriad of questions began piling up, but I pushed these aside\footnote{I later learned that it was Mike who had passed on my email, but finding Bill's email address had taken some digging. It was provided, via a colleague of Mike's at UQ and a colleague of Bill's at JCU, by the widow of a former colleague of Bill's.} and started reading the email. Now long since retired, Bill was nevertheless able to fill in the pieces of the picture I was missing. His doctoral studies had indeed started at UQ, but not at the main campus in Brisbane. Instead, he had been based in Townsville, almost a thousand miles away, at the University College of Townsville, a college of UQ. There, since 1961, at the same time as he was pursuing his doctoral studies, he had quite remarkably held a lectureship in mathematics. By 1968, he had finished writing his thesis, and was ready to graduate by 1969. Around that time, however, the college officially became James Cook University (JCU), the second university in Queensland; hence Bill was given the option to graduate with a degree from either UQ or JCU, and chose the latter. This explained why UQ did not have a copy of his thesis, as no copy was ever submitted there.

As to the matter of his supervisor, it was sporadic\footnote{Bill recalls how Baumslag sent him a letter written on a German hotel letterhead, told him he was writing from England, and had two weeks later posted the letter from South Africa.} enough that even Bill himself was not entirely certain on the matter. He thought it might have been Mike Newman initially\footnote{As mentioned earlier, Mike had been Bill's master's supervisor. The two met at a conference, and got talking because of their shared surname. Mike Newman had himself in turn been supervised by B. H. Neumann, much like Baumslag -- the two had shared a student office at the Victoria University of Manchester -- and in this case, the similarity in surnames was not coincidental: the two were first cousins once removed.}, which quickly changed when Baumslag visited Australia in 1964. Baumslag agreed to supervise Bill, but disappeared almost as soon as he had appeared. Not long thereafter, however, Baumslag sent a copy of a draft of a chapter on one-relator groups from a forthcoming book on combinatorial group theory by Magnus, Karass, and Solitar\footnote{This book, including the chapter on one-relator groups, later appeared as \cite{Magnus1966}.}. In Bill's own words, this was ``the most useful help [Baumslag] provided'', and it prompted a fruitful investigation into one-relator groups that would culminate in the Spelling Theorem.

After spending some time working with one-relator groups\footnote{Among other things, he proved a 1964 conjecture due to Baumslag, the proof of which later appeared in \cite{Newman1968b}.}, Bill was due for a sabbatical leave, and Baumslag was in 1967 able to set up a lectureship for him at Fairleigh Dickinson University in Teaneck, New Jersey. This meant that the two mathematicians were able to meet face to face again. At their very first meeting, when explaining how far he had come in proving a Spelling Theorem, Bill realised that his proof was correct. Baumslag, enthusiastically, suggested that the theorem be presented at Magnus' weekly group theory seminar at the State University in Washington Square, and this came to pass. In the audience that day were both Magnus and Solitar, two of the three who had taught Bill extensively about one-relator groups. What was more, one of these was, of course, Wilhelm Magnus, the man who had proved the \textit{Freiheitssatz}, the most significant result on one-relator groups to date. Not fazed by this, Bill presented his results, and concluded his presentation. At that point, Magnus had a remarkable reaction. He jumped to his feet, and exclaimed for the entire room to hear: ``I don't believe this! I don't believe this!''. In Bill's own words, he was saying that ``he could not believe that an unheard of mathematician, from some unknown university in outback Australia, could have come up with these results''. But the proof was correct. 

After my correspondence with Bill, there yet remained a major unresolved issue in the recovery of the thesis. I wrote an abridged email to JCU containing only the essentials of the story up to that point, for fear of the thesis slipping through my fingers in the time it would take them to read the entire story. Due to time differences, I woke up the  next morning pleasantly surprised, having discovered that I had been copied into an overnight flurry of emails sent back and forth between the archivists, heads of research, and librarians at JCU. Then, at last, the final email of the conversation stated that the thesis had been found, alive and well, and I could not have been happier. Shortly thereafter, I received an email containing the scanned thesis.

And so, a few weeks after it had begun, that chapter of the story was more or less over. I was able to read through the original proof, just as I had wanted, and I do not believe that I will ever read a proof with as much enthusiasm again. The thesis was later uploaded to JCU's online archives, and also passed on to UQ, so that they may quickly help anyone in the future digging into the same story. There was, at this point, only one small part of the story remaining. At the time of his graduation, Bill worked as a bookbinder. Indeed, at the formal opening of James Cook University in 1970 a visitors' book bound by him was present\footnote{This book lies in the Vice-Chancellor's office to this day, and remains in constant use.} and signed by Queen Elizabeth II. Due to this bookbinding interest, Bill had also decided to bind a copy of his own thesis in kangaroo leather, inlaid with a diagram of the amalgamated free product, and had it sent to Baumslag. I have been in contact with numerous people regarding the current whereabouts of this piece of paraphernalia but, as of yet, have not had any luck in tracking it down. One day, it might be found, or it may simply be lost to time. 

As for Bill Newman himself, after leaving New Jersey -- and following a six month visiting lectureship at Pahlavi University, Shiraz, Iran -- he resumed his lectureship at what had become James Cook University in 1970. He retained this position after his graduation, and would remain employed there until his retirement in 2000. A number of sabbatical years during this time are noteworthy. In 1980, he spent some months as a visiting professor at the California Institute of Technology, and worked together with Dr R. Miller at the nearby Jet Propulsion Laboratory. This interest in space travel would continue; at his next sabbatical in 1982, Bill worked in London and the Netherlands as a consultant to a company contracted by the European Space Agency. This work, on space communication and error-correcting codes, was for the 1985 interception by the \textit{Giotto} spacecraft of Halley's comet; Bill's involvement in this project subsequently made the front-page of the Townsville local newspaper. In 2000, Bill retired from his lectureship, and, at the time of the writing of this article, currently lives in the Philippines. 

\clearpage

\section{Combinatorial Group Theory}

Group presentations, the objects of interest for spelling theorems, were introduced by Walther Dyck (1856--1934), later von Dyck, in 1882 \cite{Dyck1882}, and exposed the field of group theory to a new approach for dealing with infinite groups. Dyck's approach to the subject via presentations was \textit{combinatorial} in nature, and grew naturally out of the geometric and topological approaches which permeated the topic at the time, with ideas stemming from e.g. Lie group theory \cite{Lie1880} and topological work of Klein, Fricke, and Poincar\'{e}. For example, the notion of a fundamental group of a topological space was introduced by Poincar\'{e} \cite{Poincare1895} in the mid-1890s, but to study such groups in-depth it was clear that a more combinatorial approach was called for. This proved to be the beginning of a fruitful marriage of ideas from group theory and topology.

A presentation of a group $G$ is a way of giving $G$ as a quotient of a free group. In other words, it is a way of representing elements of $G$ as equivalence classes of elements of a free group, and hence, by choosing representatives, a way to represent elements of $G$ as elements of a free group. Throughout this text, elements of free groups will always be referred to as \textit{words}, and this is no misnomer: they are words with alphabet any set of cardinality twice the rank of the free group, to accommodate formal inverse symbols. One of Dyck's many insights was that it is possible to give a presentation for \textit{any} group. This is done by specifying a number of generators and a number of defining relations. In practice, this latter quantity often turns out to be much more vital to the structure of the group than the former. If the generating alphabet is $A$ and the set of defining relations, viewed as words over the alphabet and its inverses, is $R$, then we record this presentation as $\pres{}{A}{R=1}$, or simply $\pres{}{A}{R}$, and we say that $G$ is defined by this presentation.

It is often prudent to only consider presentations with finitely many generators and defining relations. Such presentations are called finite, and a group which is defined by some finite presentation is called finitely presented, a finiteness condition on the group. It is important to note that not all groups can be finitely presented. Finite presentations allow us to deal with groups in a combinatorial way, by dealing with group elements not as abstract objects floating in a void\footnote{It may be helpful to compare this with \textit{geometric} group theory, in which group elements are studied through the way in which they act on certain geometric objects.}, but as equivalence classes of words. While in general a given element will have infinitely many words representing it in this manner, the multiplication of such words through simple concatenation is compatible with the group structure, which results in a very convenient manner of describing the group. This is especially true of finite presentations, in which case one can, for instance, describe an infinite group by only specifying a finite amount of information, a considerable improvement. This section is not meant to serve as an introduction to combinatorial group theory, and the only notion essential to understanding the remainder is that one may use formal words to represent group elements. The reader is directed to \cite{Lyndon1977} or \cite{Magnus1966} for excellent introductions to the subject, and to Stillwell \cite{Stillwell1980} for a good overview of the topological background.

Although all finite groups are finitely presented, for studying finite groups finite presentations are not always particularly helpful. Consider the example of the alternating group $A_5$, a finite group with a great deal of history in the theory of groups. It can be described uniquely in many ways; for example, it is the set of all even permutations on a set with five elements; the smallest non-abelian simple group, proved already by Galois \cite{Galois1831}; the smallest non-solvable group; or indeed the collection of rotational symmetries of the icosahedron. All of these descriptions reveal a great deal of information about the group. On the other hand, $A_5$ also admits the following presentation, which was given already in 1856 by W. R. Hamilton\footnote{The author thanks John Stillwell for bringing this fact to his attention.} \cite{Hamilton1856}, and which additionally appears in Dyck's 1882 article \cite{Dyck1882}:
\[
\pres{}{a, b}{a^5 = b^2 = (ab)^3 = 1}.
\]
Almost no information about the structure of the group is revealed directly by this presentation. In particular, it is not immediately clear that the group defined by this presentation should be finite or indeed non-trivial, let alone non-abelian simple. One can obtain similar presentations for all symmetric groups $S_n$ and alternating groups $A_n$; this was first done in 1896 by E. H. Moore\footnote{The author thanks James East for bringing this fact to his attention.} \cite{Moore1896}, and these presentations are similarly terse regarding the many properties of the groups they present. The primary target class of groups for presentations as a tool, then, lies outside of the finite\footnote{For completeness, it is true that some more modern approaches to the use of presentations to studying finite groups have been found on the computational side, as well as some specialised work by e.g. Coxeter \& Moser \cite{Coxeter-Moser}.}.

On the other hand, some infinite groups admit straightforward and illuminating presentations. For example, the free abelian group $\mathbb{Z}^2$ admits the presentation $\pres{}{a,b}{[a,b]=1}$, where $[a, b] = aba^{-1}b^{-1}$ denotes the commutator. Indeed, the group $G$ defined by this presentation is abelian (as $ab=ba$) and $2$-generated (by $a$ and $b$). Furthermore, one may check without much difficulty that no element of $G$ has finite order. It thereby follows by the classification theorem of finitely generated abelian groups that $G$ is isomorphic to the free abelian group $\mathbb{Z}^2$. Thus, in this case, the presentation captured quite succinctly all the important information about the group it presented. 

The reader may at this point feel confounded as to what the usefulness of presentations may be. One major advantage of presentations is that they may easily be manipulated, leading to new presentations for known groups. In 1908, Heinrich Tietze (1880--1964) worked from a topological viewpoint to introduce a number of transformations of presentations in \cite{Tietze1908}, which today bear his name. These permit the transformation of some given presentation into another equivalent presentation, that is, a presentation defining the same group, with one possible end goal of finding the simplest presentation for the group. For example, removing a generator which can be expressed as a word in the other generators, and subsequently replacing any occurrence of it in the presentation by this word, is a Tietze transformation. For example, the presentation $\pres{}{a,b,c}{abc = 1, cb^{-2}=1}$ defines the same group as the presentation $\pres{}{a,b}{ab(bb) = 1}$, by removing the generator $c = b^2$.  We may now see that as, in this new presentation, $a = b^{-3}$, we can remove the generator $a$, and find that our original presentation defines the same group as does $\pres{}{b}{\varnothing}$, namely the infinite cyclic group $\mathbb{Z}$.

Importantly, Tietze proved that any two presentations defining the same group can be transformed into each other through a finite sequence of applications of Tietze transformations\footnote{This proof is not constructive, and the problem of deciding whether two presentations define the same group, known as the \textit{isomorphism problem}, remains; indeed, this problem was acknowledged already in 1908 by Tietze in \cite{Tietze1908}.}. Using his transformations, Tietze \cite{Tietze1908} was able to prove that the fundamental group of a topological space is invariant under homeomorphism; this was known in part previously to Poincar\'e \cite{Poincare1895} under the additional assumption that the fundamental group be finite. Another use of presentations comes from knot theory. An important invariant of any knot is its knot group, that is, the fundamental group of the complement in $\mathbb{R}^3$ of the knot. Given any knot, it is straightforward to write down a finite presentation for this group. In 1910, Dehn \cite{Dehn1910} showed that a knot is trivial if and only if its fundamental group is abelian, hence relating the important, and very topological, problem of deciding whether a knot is trivial to solving problems about finite presentations. For more details on the topological underpinnings of presentations, the exposition given by Stillwell \cite{Stillwell1982} is an excellent starting point.

Furthermore, presentations have a key selling point: they allow us to ask questions about presentations of groups of a combinatorial nature, and, by extension, ask such questions about the groups themselves. A typical example of such a question is asking about the structure of groups given by presentations with only a single defining relation, commonly called one-relator groups. Questions of this nature can only arise in this combinatorial context, and are central to combinatorial group theory.

\section{The Word Problem}

The word problem is a central problem in the theory of presentations. Here, it will serve as a natural stepping stone from presentations to spelling theorems. Presentations let us represent group elements as words. Given a presentation and two words, it is therefore natural to ask whether the two words represent the same element. The problem of producing an effective procedure for deciding this problem in a finite time for any two input words is known as the word problem for the presentation. The word problem was introduced\footnote{More accurately, Dehn presented three problems for presentations: the \textit{identity problem}, the \textit{conjugacy problem}, and the \textit{isomorphism problem}. The identity problem is equivalent to the word problem for groups.} in 1911 by Max Dehn (1878--1952), a highly prolific student of Hilbert's, in a seminal article \cite{Dehn1911}. Coming from a topological background, he had already proved results regarding knot groups the year before, and in 1914 \cite{Dehn1914} solved the word problem for the trefoil knot group. Using an elegant approach, he was able to apply this solution of the word problem to devise a proof that the trefoil knot is not isotopic to its mirror image. This gave credence to the notion that combinatorial group theory, having grown out of topological considerations, was already becoming an incredibly useful tool for proving results on its own.

Contemporary with this work of Dehn was the Norwegian mathematician Axel Thue (1863--1922). He stated the word problem in the more general context of semigroups, and in 1914 he founded an attempted systematic approach for solving the word problem. It is important to note that even at this time there was a great deal of pessimism as to whether such an approach would ever be successful. Indeed, Dehn himself had not understated the difficulty of these decision problems in 1911, and in his own words ``solving the [word] problem for groups may be as impossible as solving all mathematical problems'' \cite{Chandler-Magnus}. He was, as we shall now see, correct in this assertion. An important branch of mathematics which appeared in the 1930s, many of whose further branches would quickly become confluent with combinatorial group theory, arose from the \textit{Entscheidungsproblem} and computability theory. One of the goals of this branch was to show the general unsolvability of certain types of decision problems. In other words, it was focused on showing that there are certain problems which cannot be solved in finite time by any effective procedure; no computer could solve them for us.

Not long thereafter, the existence of unsolvable decision problems was proved. This led to a frenzy of research on the algebraic side, in which efforts were made to encode these unsolvable decision problems in the word, conjugacy, and isomorphism problems for groups and semigroups, to prove the same problems unsolvable. This frenzy became a successful venture. In 1937, Church \cite{Church1937} published a proof of the existence of a finitely generated semigroup with unsolvable word problem; ten years later, in 1947, this was extended to the existence of a finitely presented semigroup with unsolvable word problem, independently by Post \cite{Post1947} and Markov \cite{Markov1947}. In 1950, Turing \cite{Turing1950} built on this by demonstrating the existence of a finitely presented cancellative semigroup, an intermediary structure between semigroups and groups, with unsolvable word problem. Novikov \cite{Novikov1954} proved in 1954 that the conjugacy problem is unsolvable in general for finitely presented groups. Finally, in 1955, Novikov \cite{Novikov1955} proved the existence of a finitely presented group with unsolvable word problem, placing the final nail in the coffin for the hope of a general procedure for solving the word problem for groups; it is worth mentioning that Boone \cite{Boone1958} independently reproved this result in 1958. Around this time, the isomorphism problem for finitely presented groups was also proved to be unsolvable in general by Adian \cite{Adian1955} and Rabin \cite{Rabin1958}. Thus all three of the problems posed by Dehn in 1911 had been proved unsolvable in general, and the impact they had made on the development of combinatorial group theory was undeniable.

Indeed, the word problem continued, and continues, to play a key r\^ole in combinatorial group theory in spite of these negative results. The general unsolvability of the word problem may be phrased another way: the way a word is spelled, that is, written as an element of a free group, is not in general enough to decide what kind of an element it represents in the quotient. This prompts the natural question of under which circumstances one may prove that the spelling of a word \textit{is} enough to decide what kind of an element it represents. It is answers to this question that lead us to the notion of \textit{spelling theorems}\footnote{The name \textit{spelling theorem} is originally due to B. B. Newman. As a child in the 1940s, his sometimes less than adequate spelling led to frequent canings in the headmaster's office. It was thus with a sense of delighted vindication that he decided to name his own spelling theorem.}. 

The first spelling theorem appeared in the context of \textit{surface groups}. A surface group is the fundamental group of some closed $2$-manifold (or \textit{surface}), which we for brevity will assume to be orientable. If $\mathcal{O}_g$ denotes a surface of genus $g > 0$, then the fundamental group of $\mathcal{O}_g$ admits the presentation 
\[
\pres{}{a_1, b_1, \dots, a_g, b_g}{[a_1, b_1][a_2, b_2] \cdots [a_g, b_g] = 1}
\]
where $[a_i, b_i] = a_ib_ia_i^{-1}b_i^{-1}$ as usual denotes the commutator. Importantly, the fundamental group admits a presentation with only a single defining relation. Dehn had proved already in 1911 that the word problem is decidable for surface groups using geometric methods, but in 1912 he presented what became known as \textit{Dehn's algorithm}, an algebraic method for solving the word problem for surface groups, and of which a key component is a spelling theorem. 

\begin{theorem}[Dehn's Spelling Theorem, 1912]
Let 
\[
G = \pres{}{a_1, b_1, \dots, a_g, b_g}{[a_1, b_1][a_2, b_2] \cdots [a_g, b_g] = 1}
\]
be the fundamental group of a genus $g \geq 2$ closed orientable $2$-manifold. Let $w$ be a non-trivial reduced word in the generators which represents the identity element of $G$. Then there exists a subword $u$ of $w$ such that $u$ is a subword of some cyclic conjugate of the defining relation, and such that the length of $u$ is strictly more than half the length of the defining relation.
\end{theorem}

This spelling theorem roughly states that any word representing the identity element should contain a large subword of the relator word (or some cyclic conjugate of the relator word). Through an iterative procedure, one then replaces this large subword with the shorter inverse of what remains of the relator word, and applies the spelling theorem again, this time to a shorter word. This process is guaranteed to terminate, as the word becomes shorter at every step, and the final word in the sequence is the empty word if and only if the original word represents the identity. This is enough to solve the word problem. In other words, Dehn was able to control the spelling of any word representing the identity sufficiently much to solve the word problem. This is at the heart of the notion of a spelling theorem. In honour of this, any presentation which admits a spelling theorem, once rigorously defined, is often called a Dehn presentation. The interested reader may find a lucid exposition of Dehn's algorithm in the more modern context of rewriting systems in \cite{Book1993}. We now turn from spelling theorems to the theory of one-relator groups, the final piece needed to state the B. B. Newman Spelling Theorem.

\section{One-relator Groups}

A one-relator group is one which can be presented using only a single defining relation. We have already seen examples of one-relator groups in the earlier surface groups, but there is an abundance of other examples, especially in the topological setting. For example, the knot group of the trefoil knot, or indeed any torus knot, is a one-relator group. The growth of combinatorial group theory in the mid-1900s correlates well with the development of the theory of one-relator groups. It is not far off to say that Dehn initiated this theory when he in 1928 asked his doctoral student Wilhelm Magnus (1907--1990) to prove a result about one-relator groups: the  \textit{Freiheitssatz}, or freeness theorem. Given a one-relator group, with a set of generators and a single relation, in which we may without loss of generality assume all generators appear at least once, the \textit{Freiheitssatz} says that any proper subset of the generating set generates a free group. Dehn was himself more or less convinced that this important structural result was true through geometric intuition, and presented a sketch of a geometric proof for Magnus to fill in the details of \cite{MactutorWP}. To the geometer Dehn's surprise, and perhaps slight disapproval, Magnus proved the result in 1930 using purely combinatorial methods, foregoing the geometric methods of his supervisor; indeed, Magnus recalls in \cite{Magnus1978} that, upon hearing of the algebraic nature of the proof, Dehn  exclaimed ``\textit{Da sind Sie also blind gegangen!} [So you went about it blind!]''. But the proof was correct.

What was more, Magnus was shortly thereafter able to use the \textit{Freiheitssatz} to prove that the word problem is decidable for \textit{any} one-relator group, a result which remains spectacular in its profundity to this day\footnote{Remarkably, the word problem for two-relator groups remains an open problem; this is Problem~9.26 of the Kourovka Notebook \cite{Khukhro2018}.} \cite{Magnus1932}. The theory of one-relator groups was much developed in the decades following this result, but importantly, most results were fundamentally proved using a key strategy. This strategy relies on the \textit{Magnus-Moldavanskii hierarchy}, named after Magnus, who first used it in the proof of his \textit{Freiheitssatz}, and Moldavanskii \cite{Moldavanskii1967}, who developed it into its modern form\footnote{The hierarchy was originally not studied through HNN-extensions, a notion introduced in and named after the authors of \cite{Higman1949}, but they are present in \cite{Moldavanskii1967}, and the formulation of the hierarchy through such extensions is very elegant.}. The general idea is that one may prove results about one-relator groups by induction on the length of the defining relation. Starting with any one-relator group $G_0 = \pres{}{A_0}{w_0=1}$, one associates to $G_0$ a new one-relator group $G_1 = \pres{}{A_1}{w_1 = 1}$, where $|w_1| < |w|$, i.e. $w_1$ is a shorter word than $w$ over the new finite generating set $A_1$. By a careful analysis, one finds that many properties of the group $G_1$ can now be lifted to properties of $G_0$. One may then iterate this procedure, obtaining from $G_1$ a one-relator group $G_2$, giving a sequence $G_0, G_1, G_2, \dots$ of one-relator groups. As the length $|w_i|$ is strictly decreasing, this sequence always terminates, and in fact always terminates in a one-relator group isomorphic to a direct product $F_n \times C_k$ for some finitely generated free group $F_n$ and some finite cyclic group $C_k$ (possibly trivial). As the properties of $F_n \times C_k$ are well understood, this means, modulo understanding the iteration above, that all one-relator groups can be made well understood. It is precisely in this manner that Magnus proved the word problem decidable for all one-relator groups. 

While the proofs involving the procedure above can get rather intricate, the actual process of obtaining $G_{i+1}$ from $G_i$ is very straightforward to write down. We give an example of this. Suppose we start with the one-relator group $G_0$ given to us as $G_0 = \pres{}{a,b}{abab^{-1}a^{-1}ba^{-1}b^{-1} = 1}$, i.e. we have 
\[
w_0 = abab^{-1}a^{-1}ba^{-1}b^{-1}.
\]
We notice that the \textit{exponent sum} $\sigma_b(w_0)$ of $b$ in $w_0$, i.e. the total sum of the exponents of the occurrences of $b$ in $w_0$, is $\sigma_b(w_0) = 0$. Using this, we first obtain a word $w_0'$ from $w_0$ by replacing every occurrence of the letter $a$ (resp. $a^{-1}$) by $a_i$ (resp. $a_i^{-1}$), where $i$ is the sum of the exponents of the $b$'s appearing before that occurrence in $w_0$. In our example, the reader may verify that we obtain the word 
\[
w_0' = a_0 ba_1 b^{-1}a_{0}^{-1}ba_1^{-1}b^{-1}.
\]
We now remove from $w_0'$ all occurrences of our exponent-sum zero letter $b$, leaving us with the word
\[
w_1 = a_0 a_1 a_0^{-1} a_{1}^{-1}.
\]
The set of generators $A_1$ is now just the set of letters $a_i$ which appear in this word, which in this case is $A_1 = \{ a_0, a_1 \}$, and we find 
\[
G_1 = \pres{}{a_0, a_1}{a_0 a_1 a_0^{-1} a_{1}^{-1} = 1} = \pres{}{a_0, a_1}{[a_0, a_1] = 1}.
\]
Thus $G_1$ is isomorphic to the free abelian group $\mathbb{Z}^2$. For ease of notation, we make a change of generating set, and write $G_1 = \pres{}{c,d}{cdc^{-1}d^{-1} = 1}$. We now repeat the same process -- this time using, say, the letter $d$, which has exponent sum zero in $cdc^{-1}d^{-1}$ -- and obtain the group 
\[
G_2 = \pres{}{c_0, c_1}{c_0 c_1^{-1} = 1}
\]
as the reader may verify. Using a Tietze transformation, we find that this group is isomorphic to the free group on one generator, i.e. the infinite cyclic group $\mathbb{Z}$, completing the procedure. 

Of course, in a general one-relator group $\pres{}{A}{w=1}$, there may not be any letter $a \in A$ such that the exponent sum $\sigma_a(w)$ is $0$. In this case, one may without much difficulty modify the defining relation slightly, by the introduction of a new generator, in a way which simultaneously ensures the group-theoretic properties are virtually unchanged as well as ensures that the new defining relation contains a generator with exponent sum $0$. Furthermore, we can generalise the above from two generators to arbitrarily many in the obvious way, i.e. by choosing a single generator with exponent sum zero, and re-indexing the other generators as above.  Hence, this allows the above procedure to apply to all one-relator groups. For a full and rigorous treatment of the hierarchy (with more examples), we refer the reader to the excellent exposition found in McCool \& Schupp \cite{McCool1973}.

One key player which entered the world of one-relator groups around 1960, and made extensive use of the Magnus-Moldavanskii hierarchy, was Gilbert Baumslag (1933--2014), a student of B. H. Neumann's. His work proved invaluable to developing the field into its modern form, and we will now focus on one central aspect of this development. There is a natural split of the theory of one-relator groups into two areas of study: the \textit{torsion} and \textit{torsion-free} case. A group has torsion if it has non-trivial elements of finite order, and is torsion-free otherwise. In 1960, a simple combinatorial characterisation on the defining relation of a given one-relator group was found, which determines whether the group has torsion or not: a one-relator group $\pres{}{A}{w = 1}$ has torsion if and only if $w = r^n$ for some word $r$ and $n > 1$, and any element of finite order is conjugate to some power of $r$ \cite{Karrass1960}. In other words, the only torsion is the obvious one. At this point in time, the general consensus was that the torsion case seemed significantly more well-behaved and easier to study than its torsion-free counterpart. This prompted Baumslag to pose a series of conjectures regarding one-relator groups with torsion.

\begin{conjecture}[Baumslag's Conjectures]
Let $G = \pres{}{A}{w^n = 1}$ be a finitely generated one-relator group with torsion. Then \begin{enumerate}
\item $G$ is residually finite \cite{Baumslag1964};
\item $G$ is virtually free-by-cyclic \cite{Baumslag1967};
\item $G$ is coherent \cite{Baumslag1974}.
\end{enumerate}
\end{conjecture}

Here, a \textit{residually finite} group is one in which the intersection of all subgroups of finite index is trivial; a \textit{virtually free-by-cyclic} group is one with a finite index subgroup which is an extension of a free group by a cyclic group; and a \textit{coherent} group is one in which every finitely generated subgroup is finitely presented. All of these notions are rather specialised, and will not be expanded upon further in the course of this article. However, it bears mentioning that all such properties are very strong; for example, none of the properties hold for two-relator groups, in general, and there exist non-residually finite and not virtually free-by-cyclic torsion-free one-relator groups. One of the major reasons that such strong conjectures can be made in the torsion case is directly because of the B. B. Newman Spelling Theorem. As we shall see, it is in fact one of the foundational results behind the recent affirmative resolutions to the first and the third of the three above conjectures.

\section{The B. B. Newman Spelling Theorem}

With a sufficient \textit{aper\c{c}u} of some of the relevant theory of one-relator groups provided, it is now possible to present the namesake of this article: a spelling theorem, much like Dehn's spelling theorem for surface groups.

\begin{theorem}[The B. B. Newman Spelling Theorem, 1968]
Let $G = \pres{}{A}{R^n = 1}$ be a one-relator group with torsion such that $R$ is cyclically reduced and not a proper power. Let $w$ be a non-trivial reduced word representing the identity element of $G$. Then $w$ contains a subword $u$ such that either $u$ or $u^{-1}$ is a subword of $R^n$, and such that the length of $u$ is strictly more than $\frac{n-1}{n}$ times the length of $R^n$.  
\end{theorem}

If the reader recalls the formal statement of Dehn's spelling theorem for surface groups, they will pleasantly note many similarities. The theorem first made an appearance in a 1968 bulletin article \cite{Newman1968}\footnote{The Spelling Theorem is written in a slightly different form in the bulletin article and Newman's thesis. In that form, the theorem instead gives similar information as above about the spelling of non-trivial words equal to each other.}, but can be found in almost any context in which one-relator groups are discussed. In particular, a very accessible proof of the theorem can be found in either of \cite{Lyndon1977, McCool1973}, and this proof takes full advantage of the Magnus-Moldavanskii hierarchy. We give an idea of the proof.

The proof, as with most proofs of theorems about one-relator groups, is by induction on the length of the relator word. One may first observe some base cases. If $G = \pres{}{A}{R^n=1}$ only has a single generator, then $G = \pres{}{a}{a^n = 1}$ is a finite cyclic group. A non-trivial word $a^i$ ($i \in \mathbb{Z} \setminus \{ 0 \}$) in the free group on the single generator $a$ is equal to the identity element of $G$ if and only if $i \equiv 0 \mod n$. Thus, any such word $a^i$ contains either $a^n$ or $a^{-n}$ as a subword, and the statement of the theorem is quickly checked to hold. If instead $G = \pres{}{A}{R^n =1}$ has more than a single generator, and at least two generators appear in $R$, then one may apply the Magnus-Moldavanskii procedure illustrated earlier. At each step of shortening the relation, one must now ensure that one maintains control over words equal to the identity; originally, this was done via a clever shifting of the subscripts of the new generators obtained at each stage of the hierarchy, at which point the inductive hypothesis may be applied. The original proof in \cite{Newman1968b} is not long, taking up six typewritten pages, and the main issue to be resolved are the many different cases of potential exponent sums; the proof given in \cite{Lyndon1977} is very similar.

There are many important consequences of the Spelling Theorem. The first is that it allows for Dehn's algorithm to be used to solve the word problem for one-relator groups with torsion, just as it was used to solve the word problem for surface groups. In view of Magnus' 1932 solution of the word problem for \textit{all} one-relator groups, this consequence is perhaps not particularly exciting\footnote{However, the Spelling Theorem does demonstrate that the word problem for one-relator groups with torsion is decidable in \textit{linear time}, whereas the exact time complexity of Magnus' solution remains unknown in general.}. A stronger consequence, already observed in \cite{Newman1968}, is that the conjugacy problem, one of Dehn's three original problems, is decidable for one-relator groups with torsion. Furthermore, Pride \cite{Pride1977b} relies heavily on the Spelling Theorem in his celebrated article showing that the isomorphism problem is solvable for two-generated one-relator groups with torsion, and therein presents an entirely combinatorial argument of this fact, free of geometry. The well-behaved combinatorial nature of the setting has made the theorem straightforward to generalise to much broader geometric situations \cite{Howie1984}. 

The Spelling Theorem can also be directly applied to prove that one-relator groups with torsion are \textit{hyperbolic}, a property introduced by Gromov \cite{Gromov1987}, and the conjugacy problem is decidable for any hyperbolic group. This hyperbolicity plays a crucial r\^ole in Wise's \cite{Wise2009} highly geometric resolution of Baumslag's conjecture on the residual finiteness of one-relator groups with torsion, where it is shown that any such group admits a quasiconvex hierarchy; in combination with hyperbolicity, this is sufficient to guarantee virtual specialness, a property more than sufficient for residual finiteness. The interested reader may consult \cite{Wise2012} for a full exposition of the first resolution. The hyperbolicity plays a similar r\^ole in the recent proof due to Louder \& Wilton \cite{Louder2018} (and, independently, Wise\footnote{The author thanks Daniel Wise for bringing this fact to his attention.} \cite{Wise2018}) that all one-relator groups with torsion are coherent, thus also resolving the third of Baumslag's conjectures. It is therefore clear that the Spelling Theorem was an important part of resolving both the first and the third of Baumslag's three conjectures in the affirmative. Baumslag's second conjecture remains open. 

Hence, serving as a stepping stone from combinatorial to geometric methods, as well as an important tool for one-relator groups in distinguishing the torsion case from the torsion-free, it is clear that the Spelling Theorem remains one of the most important results proved in combinatorial group theory to date, and continues to play a key and active r\^ole in modern research on the topic.

\section*{Acknowledgements}
The author thanks Mike Newman and Christine Dalais for all their assistance in uncovering the missing thesis, and Asaf Karagila for encouraging the write-up of this material. The author also gratefully acknowledges the anonymous referee for their many detailed and helpful comments. Finally, the author wishes to give a special thanks to Bill Newman for the information and riveting anecdotes he has so kindly contributed, without which this essay could not have been written.

\bibliography{BBNewmanRevised} 

\providecommand{\bysame}{\leavevmode\hbox to3em{\hrulefill}\thinspace}
\providecommand{\MR}{\relax\ifhmode\unskip\space\fi MR }
\providecommand{\MRhref}[2]{%
  \href{http://www.ams.org/mathscinet-getitem?mr=#1}{#2}
}
\providecommand{\href}[2]{#2}
\begin{thebibliography}{New68b}

\bibitem[Adi55]{Adian1955}
S.~I. Adian, \emph{Algorithmic unsolvability of problems of recognition of
  certain properties of groups}, Dokl. Akad. Nauk SSSR (N.S.) \textbf{103}
  (1955), 533--535. \MR{0081851}

\bibitem[Bau64]{Baumslag1964}
Gilbert Baumslag, \emph{Groups with one defining relator}, J. Austral. Math.
  Soc. \textbf{4} (1964), 385--392. \MR{0172909}

\bibitem[Bau67]{Baumslag1967}
\bysame, \emph{Residually finite one-relator groups}, Bull. Amer. Math. Soc.
  \textbf{73} (1967), 618--620. \MR{212078}

\bibitem[Bau74]{Baumslag1974}
\bysame, \emph{Some problems on one-relator groups}, Proceedings of the
  {S}econd {I}nternational {C}onference on the {T}heory of {G}roups
  ({A}ustralian {N}at. {U}niv., {C}anberra, 1973), 1974, pp.~75--81. Lecture
  Notes in Math., Vol. 372. \MR{0364463}

\bibitem[BO93]{Book1993}
Ronald~V. Book and Friedrich Otto, \emph{String-rewriting systems}, Texts and
  Monographs in Computer Science, Springer-Verlag, New York, 1993. \MR{1215932}

\bibitem[Boo58]{Boone1958}
William~W. Boone, \emph{The word problem}, Proc. Nat. Acad. Sci. U.S.A.
  \textbf{44} (1958), 1061--1065. \MR{101267}

\bibitem[Chu37]{Church1937}
Alonzo Church, \emph{Combinatory logic as a semi-group}, Bulletin of the
  American Mathematical Society \textbf{43} (1937), 333.

\bibitem[CM80]{Coxeter-Moser}
H.~S.~M. Coxeter and W.~O.~J. Moser, \emph{Generators and relations for
  discrete groups}, fourth ed., Ergebnisse der Mathematik und ihrer
  Grenzgebiete [Results in Mathematics and Related Areas], vol.~14,
  Springer-Verlag, Berlin-New York, 1980. \MR{562913}

\bibitem[CM82]{Chandler-Magnus}
Bruce Chandler and Wilhelm Magnus, \emph{The history of combinatorial group
  theory}, Studies in the History of Mathematics and Physical Sciences, vol.~9,
  Springer-Verlag, New York, 1982, A case study in the history of ideas.
  \MR{680777}

\bibitem[Deh10]{Dehn1910}
M.~Dehn, \emph{\"{U}ber die {T}opologie des dreidimensionalen {R}aumes}, Math.
  Ann. \textbf{69} (1910), no.~1, 137--168. \MR{1511580}

\bibitem[Deh11]{Dehn1911}
\bysame, \emph{\"{U}ber unendliche diskontinuierliche {G}ruppen}, Math. Ann.
  \textbf{71} (1911), no.~1, 116--144. \MR{1511645}

\bibitem[Deh14]{Dehn1914}
\bysame, \emph{Die beiden {K}leeblattschlingen}, Math. Ann. \textbf{75} (1914),
  no.~3, 402--413.

\bibitem[Dyc82]{Dyck1882}
Walther Dyck, \emph{Gruppentheoretische {S}tudien}, Math. Ann. \textbf{20}
  (1882), no.~1, 1--44.

\bibitem[Gal31]{Galois1831}
\'E. Galois, \emph{M\'emoire sur les conditions de r\'esolubilit\'e des
  \'equations par radicaux}, Proc\`es-Verbaux de l’Acad\'emie des Sciences de
  l'Institut de France \textbf{IX} (1828--31).

\bibitem[Gro87]{Gromov1987}
M.~Gromov, \emph{Hyperbolic groups}, Essays in group theory, Math. Sci. Res.
  Inst. Publ., vol.~8, Springer, New York, 1987, pp.~75--263. \MR{919829}

\bibitem[Ham56]{Hamilton1856}
Sir William~Rowan Hamilton, \emph{Memorandum respecting a new system of roots
  of unity}, The London, Edinburgh, and Dublin Philosophical Magazine and
  Journal of Science \textbf{12} (1856), no.~81, 446--446.

\bibitem[HNN49]{Higman1949}
Graham Higman, B.~H. Neumann, and Hanna Neumann, \emph{Embedding theorems for
  groups}, J. London Math. Soc. \textbf{24} (1949), 247--254. \MR{32641}

\bibitem[HP84]{Howie1984}
J.~Howie and S.~J. Pride, \emph{A spelling theorem for staggered generalized
  {$2$}-complexes, with applications}, Invent. Math. \textbf{76} (1984), no.~1,
  55--74. \MR{739624}

\bibitem[KM18]{Khukhro2018}
E.~I. Khukhro and V.~D. Mazurov (eds.), \emph{The {K}ourovka notebook}, Sobolev
  Institute of Mathematics. Russian Academy of Sciences. Siberian Branch,
  Novosibirsk, 2018, Unsolved problems in group theory, Nineteenth edition [
  MR0204500], March 2019 update. \MR{3981599}

\bibitem[KMS60]{Karrass1960}
A.~Karrass, W.~Magnus, and D.~Solitar, \emph{Elements of finite order in groups
  with a single defining relation}, Comm. Pure Appl. Math. \textbf{13} (1960),
  57--66. \MR{124384}

\bibitem[Lie80]{Lie1880}
Sophus Lie, \emph{Theorie der {T}ransformationsgruppen {I}}, Math. Ann.
  \textbf{16} (1880), no.~4, 441--528. \MR{1510035}

\bibitem[LS77]{Lyndon1977}
Roger~C. Lyndon and Paul~E. Schupp, \emph{Combinatorial group theory},
  Springer-Verlag, Berlin-New York, 1977, Ergebnisse der Mathematik und ihrer
  Grenzgebiete, Band 89.

\bibitem[LW18]{Louder2018}
Larsen Louder and Henry Wilton, \emph{One-relator groups with torsion are
  coherent}, Math. Res. Lett., to appear (2018).

\bibitem[Mag32]{Magnus1932}
W.~Magnus, \emph{Das {I}dentit\"{a}tsproblem f\"{u}r {G}ruppen mit einer
  definierenden {R}elation}, Math. Ann. \textbf{106} (1932), no.~1, 295--307.
  \MR{1512760}

\bibitem[Mag79]{Magnus1978}
Wilhelm Magnus, \emph{Max {D}ehn}, Math. Intelligencer \textbf{1} (1978/79),
  no.~3, 132--143. \MR{505030}

\bibitem[Mar47]{Markov1947}
A.~Markov, \emph{The impossibility of certain algorithms in the theory of
  associative systems. {II}}, Doklady Akad. Nauk SSSR (N.S.) \textbf{58}
  (1947), 353--356. \MR{0023208}

\bibitem[MKS66]{Magnus1966}
Wilhelm Magnus, Abraham Karrass, and Donald Solitar, \emph{Combinatorial group
  theory: {P}resentations of groups in terms of generators and relations},
  Interscience Publishers [John Wiley \& Sons, Inc.], New York-London-Sydney,
  1966. \MR{0207802}

\bibitem[Mol67]{Moldavanskii1967}
D.~I. Moldavanski\u{\i}, \emph{Certain subgroups of groups with one defining
  relation}, Sibirsk. Mat. \v{Z}. \textbf{8} (1967), 1370--1384. \MR{0220810}

\bibitem[Moo96]{Moore1896}
E.~H. Moore, \emph{Concerning the abstract groups of order k! and 1/2 k!
  holohedrically isomorphic with the symmetric and the alternating
  substitution-groups on k letters}, Proceedings of the London Mathematical
  Society \textbf{s1-28} (1896), no.~1, 357--367.

\bibitem[MS73]{McCool1973}
James McCool and Paul~E. Schupp, \emph{On one relator groups and {${\rm HNN}$}
  extensions}, J. Austral. Math. Soc. \textbf{16} (1973), 249--256, Collection
  of articles dedicated to the memory of Hanna Neumann, II. \MR{0338186}

\bibitem[New68a]{Newman1968b}
B.~B. Newman, \emph{Some aspects of one-relator groups}, Ph.D. thesis,
  University College of Townsville, 1968.

\bibitem[New68b]{Newman1968}
\bysame, \emph{Some results on one-relator groups}, Bull. Amer. Math. Soc.
  \textbf{74} (1968), 568--571.

\bibitem[Nov54]{Novikov1954}
P.~S. Novikov, \emph{Unsolvability of the conjugacy problem in the theory of
  groups}, Izv. Akad. Nauk SSSR. Ser. Mat. \textbf{18} (1954), 485--524.
  \MR{0075196}

\bibitem[Nov55]{Novikov1955}
\bysame, \emph{Ob algoritmi\v{c}esko\u{\i} nerazre\v{s}imosti problemy
  to\v{z}destva slov v teorii grupp}, Trudy Mat. Inst. im. Steklov. no. 44,
  Izdat. Akad. Nauk SSSR, Moscow, 1955.

\bibitem[OR08]{MactutorWP}
J.~J. O'Connor and E.~F. Robertson, \emph{Word problems}, MacTutor History of
  Mathematics (2008).

\bibitem[Poi10]{Poincare1895}
Henri Poincar\'{e}, \emph{Papers on topology}, History of Mathematics, vol.~37,
  American Mathematical Society, Providence, RI; London Mathematical Society,
  London, 2010, {{\i}t Analysis situs} and its five supplements, Translated and
  with an introduction by John Stillwell.

\bibitem[Pos47]{Post1947}
Emil~L. Post, \emph{Recursive unsolvability of a problem of {T}hue}, J.
  Symbolic Logic \textbf{12} (1947), 1--11. \MR{20527}

\bibitem[Pri77]{Pride1977b}
Stephen~J. Pride, \emph{The isomorphism problem for two-generator one-relator
  groups with torsion is solvable}, Trans. Amer. Math. Soc. \textbf{227}
  (1977), 109--139. \MR{430085}

\bibitem[Rab58]{Rabin1958}
Michael~O. Rabin, \emph{Recursive unsolvability of group theoretic problems},
  Ann. of Math. (2) \textbf{67} (1958), 172--194. \MR{110743}

\bibitem[Sti80]{Stillwell1980}
John~C. Stillwell, \emph{Classical topology and combinatorial group theory},
  Graduate Texts in Mathematics, vol.~72, Springer-Verlag, New York-Berlin,
  1980. \MR{602149}

\bibitem[Sti82]{Stillwell1982}
John Stillwell, \emph{The word problem and the isomorphism problem for groups},
  Bull. Amer. Math. Soc. (N.S.) \textbf{6} (1982), no.~1, 33--56. \MR{634433}

\bibitem[Tie08]{Tietze1908}
Heinrich Tietze, \emph{\"{U}ber die topologischen {I}nvarianten
  mehrdimensionaler {M}annigfaltigkeiten}, Monatsh. Math. Phys. \textbf{19}
  (1908), no.~1, 1--118. \MR{1547755}

\bibitem[Tur50]{Turing1950}
A.~M. Turing, \emph{The word problem in semi-groups with cancellation}, Ann. of
  Math. (2) \textbf{52} (1950), 491--505. \MR{37294}

\bibitem[Wis09]{Wise2009}
Daniel~T. Wise, \emph{Research announcement: the structure of groups with a
  quasiconvex hierarchy}, Electron. Res. Announc. Math. Sci. \textbf{16}
  (2009), 44--55. \MR{2558631}

\bibitem[Wis12]{Wise2012}
\bysame, \emph{From riches to raags: 3-manifolds, right-angled {A}rtin groups,
  and cubical geometry}, CBMS Regional Conference Series in Mathematics, vol.
  117, Published for the Conference Board of the Mathematical Sciences,
  Washington, DC; by the American Mathematical Society, Providence, RI, 2012.
  \MR{2986461}

\bibitem[Wis18]{Wise2018}
\bysame, \emph{Coherence, local-indicability and non-positive immersions},
  submitted (2018).

\end{thebibliography}
\bibliographystyle{amsalpha}

\end{document}